\documentclass{amsart}

\usepackage{amssymb,amsmath,amsthm,color}

\usepackage[bookmarksopen=true,bookmarksnumbered=true, bookmarksopenlevel=2, colorlinks=true,linkcolor=mblue,
citecolor=mblue,urlcolor=mblue, pdfstartview=FitH]{hyperref}

\definecolor{mblue}{rgb}{0,0,.8}

\newcommand{\N}{\mathbb N}
\newcommand{\Z}{\mathbb Z}
\newcommand{\Q}{\mathbb Q}

\newcommand{\A}{\mathfrak A}
\newcommand{\B}{\mathfrak B}

\newcommand{\uni}{\mathcal U}
\newcommand{\known}{\sim_+}
\newcommand{\notknown}{\sim_-}

\newtheorem{thm}{Theorem}

 \DeclareMathOperator{\Nm}{N} 

\def\dash---{\thinspace---\hskip.16667em\relax}

\begin{document}

\title[New invariants for  shift equivalence]{On some new invariants for shift equivalence for shifts of finite type.}
\author{S\o ren Eilers, Ian Kiming}
\address{Department of Mathematics, University of Copenhagen, Universitetsparken 5, DK-2100 Copenhagen \O ,
Denmark.} \email{\href{mailto:eilers@math.ku.dk}{eilers@math.ku.dk}}
\email{\href{mailto:kiming@math.ku.dk}{kiming@math.ku.dk}}

\begin{abstract} We introduce a new computable invariant for strong shift equivalence of shifts of finite type. The invariant is based on an invariant introduced by Trow, Boyle, and Marcus, but has the advantage of being readily computable.

We summarize briefly a large-scale numerical experiment aimed at deciding strong shift equivalence for shifts of finite type given by irreducible $2\times 2$-matrices with entry sum less than 25, and give examples illustrating to power of the new invariant, i.e., examples where the new invariant can disprove strong shift equivalence whereas the other invariants that we use can not.
\end{abstract}

\maketitle

\section{Introduction.}\label{intro} The shifts of finite type (SFTs) form an important class of symbolic dynamical systems which has fundamental applications in mathematics, physics and computer science. The classification problem for irreducible SFTs up to conjugacy is generally believed to be undecidable as indicated by the examples of Kim and Roush \cite{kim-roush_williams} demonstrating the difference between Williams' concepts \cite{williams} of \emph{shift equivalence} and \emph{strong shift equivalence}. Indeed, shift equivalence is decidable \cite{kim-roush_decidable}, but it is the more elusive strong shift equivalence which encodes this significant problem, and one can no longer hope that these properties are one and the same. Furthermore, the procedure in \cite{kim-roush_decidable} is not readily implementable by computer algebra systems, and hence, unfortunately, of little practical use when trying to determine what strong shift equivalences exist on a large body of matrices.
\smallskip

When attempting to prove that two matrices \textbf{fail} to be strongly shift equivalent one has access to a large and very diverse family of invariants developed over the last decades, see \cite{lind-marcus} for a summary. Most of these invariants are efficiently computable and comparable as they take the form of algebraic numbers, finitely generated groups etc. An important invariant was developed by Trow in \cite{trow} and Boyle, Marcus, Trow in \cite{bmt}. This invariant takes the form of the class of a certain ideal in a certain integral domain (see below in section \ref{inv} for details). We shall refer to it as the BMT invariant for brevity. The fact that the BMT invariant is not readily computable was the starting point and prime motivation behind the present paper. \smallskip

With this invariant as basis we introduce here $2$ new invariants: The first is defined and computable under a slight technical restriction and is proved to be equivalent to the BMT invariant. The second invariant is defined unconditionally, is possibly weaker than the BMT invariant, but has the advantage of being computable.
\smallskip

More precisely, the BMT invariant takes the form of the class (defined in the usual way) of a certain ideal of the ring $\Z[1/\lambda]$ where $\lambda$ is a certain algebraic integer. Under a technical restriction involving the conductor of the order $O:=\Z[\lambda]$ our first new invariant is defined as a certain element of the Picard group of $O$. This new invariant is shown to be equivalent to the BMT invariant when it is defined. Our second invariant is defined unconditionally as a certain element in the class group of the algebraic number field $\Q(\lambda)$. This second invariant is weaker than the BMT invariant in the sense that equality of BMT invariants implies equality of our second invariants, but it is computable by standard algorithms in algebraic number theory as implemented for instance in the computer algebra package Magma, \cite{magma}. This leads to an algorithmic approach to testing this necessary condition which, as we shall see, is quite efficient in disproving strong shift equivalence where all other invariants fail.
\smallskip

We have combined this contribution with other tools that are already described in the literature to perform a complete analysis of the question of strong shift equivalence in a limited universe of SFTs, given by irreducible integer valued $2\times 2$ matrices with an entry sum less than or equal to $25$. Building on a project by O.\ Lund Jensen \cite{lund_jensen} and using standard database tools we have recorded invariants for all matrices in this universe (with
the purpose of telling isomorphism classes SFTs apart) and concrete strong shift equivalences (with the purpose of identifying isomorphism classes).
\smallskip

The net result of these efforts can be summarized as follows. There are $17250$ matrices in the universe described, and hence $148772625$ potential questions of the type `are matrices $A$ and $B$ strong shift equivalent?' We can answer $99.99$ \% of these questions by this approach. We will briefly summarize the methods and results of this project in section \ref{comp} below.

\section{New invariants.}\label{inv} Let $S$ be an $n\times n$ matrix with non-negative, integral coefficients. We call $S$ irreducible if for every $(i,j)$ there is $k\ge 0$ such that the $(i,j)$'th entry of $S^k$ is positive ($S^0$ is defined to be the identity matrix; in other words, the irreducibility condition is empty for the diagonal entries of $S$). Irreducibility of $S$ corresponds to irreducibility of the associated SFT in the sense that any pair of legal words $u,v$ can be interpolated by a third word, $w$, to obtain a legal word $uwv$.
\smallskip

Under these conditions, one knows, cf.\ for instance Theorem 4.2.3 (the Perron--Frobenius Theorem) of \cite{lind-marcus}, that $S$ has a positive eigenvalue $\lambda$ that occurs with multiplicity $1$ in the characteristic polynomial of $S$ and whose corresponding eigenspace is $1$-dimensional, and is such that $|\mu| \le \lambda$ for any other eigenvalue $\mu$. Further, the eigenspace corresponding to $\lambda$ is generated by a positive eigenvector, i.e., a vector with positive coordinates. This uniquely determined eigenvalue is referred to as the Perron eigenvalue of $S$.
\smallskip

If now $\lambda$ is the Perron eigenvalue of $S$ there is a corresponding eigenvector
$$
{\mathbf v} = (v_1,\ldots,v_n)
$$
with coordinates $v_1,\ldots,v_n$ in the ring $\Z[\lambda]$. As the eigenspace corresponding to $\lambda$ is $1$-dimensional, the vector ${\mathbf v}$ is uniquely determined up to multiplication with a non-zero element of the algebraic number field $\Q(\lambda)$.

We can then define the BMT invariant of $S$ as the class $\mathcal{I}(S)$ of the ideal:
$$
R v_1 + \ldots + R v_n
$$
of the ring $R:= \Z[1/\lambda]$. Here, `class' has the usual meaning: Ideals $C$ and $D$ of $R$ are called equivalent if there is $\xi\in \Q(\lambda)$ such that $\xi C = D$.
\smallskip

The BMT invariant is an invariant because of the following statement that follows from Theorem 12.1 in section 12.3 of \cite{lind-marcus} (see also Theorem 6.1 of \cite{bmt}, as well as \cite{trow}): Suppose that $S$ and $T$ are matrices with integral, non-negative entries that are irreducible in the above sense and have the same Perron eigenvalue. Then if the SFTs attached to $S$ and $T$ are strongly shift equivalent, we have $\mathcal{I}(S) = \mathcal{I}(T)$.
\smallskip

In Theorem \ref{t1} below we introduce $2$ new invariants. The first of these is not always defined, but when it is defined for both $S$ and $T$ it coincides for $S$ and $T$ if and only if $\mathcal{I}(S) = \mathcal{I}(T)$, and is in this sense equivalent to the BMT invariant. This invariant takes values in the Picard group of the ring $\Z[\lambda]$.

The second invariant is always defined, takes values in the class group of the algebraic number field $\Q(\lambda)$, and is weaker than the BMT invariant in the sense that $\mathcal{I}(S) = \mathcal{I}(T)$ implies that the second invariants of $S$ and $T$ coincide.
\medskip

Now let us begin to define these new invariants and prepare Theorem \ref{t1} below. We will work with the following slightly more general setup and notation:
\medskip

\begin{tabular}{lll} $K$ &:& an algebraic number field, i.e., a finite extension of $\Q$\\
$O_K$ &:& the ring of algebraic integers in $K$\\
$\mathrm{Cl}(O_K)$ &:& the class group of $O_K$\\
$\lambda$ &:& a element of $O_K$ of the same degree over $\Q$ as $K$
\end{tabular}
\medskip

In other words, the assumption on $\lambda$ is that $K=\Q(\lambda)$. Then the ring
$$
O := \Z[\lambda]
$$
is of finite index in $O_K$, i.e., is an order of $K$. As a reference for the general theory of orders in algebraic number fields, see for instance Chap.\ 1, \S 12 of \cite{neukirch}.

In particular, attached to the ring $O$ is a certain ideal of $O$, the so-called conductor of $O$. Ideals of $O$ prime to the conductor have unique factorizations into products of prime ideals. Attached to $O$ is the Picard group $\mathrm{Pic}(O)$ of invertible ideals modulo principal ones; the Picard group coincides with the class group $\mathrm{Cl}(O_K)$ if $O=O_K$. The class group $\mathrm{Cl}(O_K)$ is a canonical quotient of $\mathrm{Pic}(O)$.

If $C$ is an ideal of $O_K$ we shall denote by $[C]$ the class of $C$ in $\mathrm{Cl}(O_K)$; similarly, if $C$ is an invertible ideal of $O$, the symbol $[C]$ denotes the class of $C$ in $\mathrm{Pic}(O)$.
\medskip

We consider the ring $\Z[1/\lambda]$. This ring is in fact the localization $O_{(M)}$ of $O$ with respect to the multiplicatively closed system
$$
M := \{ 1,\lambda,\lambda^2,\lambda^3,\ldots \} ~.
$$

The claim follows immediately once we notice that $\lambda \in \Z[1/\lambda]$: For $\lambda$ satisfies a polynomial equation:
$$
\lambda^n + a_1\lambda^{n-1} + \ldots + a_n = 0
$$
where $n:=[K:\Q]$ and the $a_i$ are integers. It follows that:
$$
\lambda = -a_1 - \ldots - a_n\cdot\frac{1}{\lambda^{n-1}} \in \Z[1/\lambda] ~.
$$

For ideals in any one of the rings we are considering above, we have the usual equivalence relation denoted by $\sim$, and defined by: $C\sim D$ if and only if there exists $\xi\in K^{\times}$ such that $\xi C = D$.
\medskip

We will now for the remainder of this section assume that we are given two ideals $\A$ and $\B$ of the ring $O_{(M)}$. Since $O_{(M)}$ is a localization of $O$ we know by general theory, see \cite{zs}, Chap. IV, \S 10, p.\ 223, that every ideal of $O_{(M)}$ is extended from an ideal of $O$, in other words, that there are ideals $A$ and $B$ of $O$ such that:
$$
\A = O_{(M)}\cdot A ~,\qquad \B = O_{(M)}\cdot B ~.
$$

We fix such ideals $A$ and $B$.
\smallskip

If $C$ and $D$ are ideals of $O$ we employ the usual notation $(C : D)$ to denote the fractional ideal:
$$
(C : D) := \{ \xi\in K \mid ~ \xi D \subseteq C \} ~.
$$

\begin{thm}\label{t1} (i). Retaining the above notation, we have $\A\sim \B$ if and only if there are elements $x\in (A:B)$ and $y\in (B :A)$ such that:
$$
xy = \lambda^k
$$
for some non-negative integer $k$.
\medskip

\noindent (ii). Suppose that the ideals $A$, $B$, and $O\cdot\lambda$ are all prime to the conductor of $O$ and let $\{ P_i\}_{i\in I}$ be the set of prime divisors in $O$ of the ideal $O\cdot\lambda$.
\smallskip

Then $\A\sim \B$ if and only if:
$$
[A] \equiv [B] \mod{\langle [P_i] \mid ~i\in I\rangle}
$$
in $\mathrm{Pic}(O)$.
\medskip

\noindent (iii). In any case, if $\{ Q_j\}_{j\in J}$ denotes the set of prime divisors in $O_K$ of the ideal $O_K\cdot\lambda$ then a necessary condition for $\A\sim \B$ is that:
$$
[O_K\cdot A] \equiv [O_K\cdot B] \mod{\langle [Q_j] \mid ~j\in J\rangle}
$$
in $\mathrm{Cl}(O_K)$.
\end{thm}
\begin{proof} We first prove the `only if' parts of {\it (i)} and {\it (ii)} as well as part {\it (iii)} simultaneously. So, suppose that $\A\sim \B$. Since $K$ is the field of fractions of $O$ there are nonzero elements $\alpha,\beta\in O$ such that $\alpha \A = \beta \B$, i.e., such that:
$$
O_{(M)} \cdot \alpha A = O_{(M)} \cdot \beta B ~.
$$

Now, $O$ is a Noetherian ring so the ideals $A$ and $B$ are finitely generated $O$-modules. Write:
$$
A = \sum_{\sigma\in S} O\cdot a_{\sigma} ~,\qquad B = \sum_{\tau\in T} O\cdot b_{\tau}
$$
with certain elements $a_{\sigma}\in A$, $b_{\tau}\in B$, and finite index sets $S$ and $T$.
\smallskip

For each $\sigma\in S$ we then have $\alpha a_{\sigma} \in O_{(M)}\cdot \beta B$ and so there is an element $s_{\sigma}\in M$ such that:
$$
s_{\sigma} \alpha a_{\sigma} \in \beta B ~.
$$

Similarly, there is for each $\tau\in T$ an element $t_{\tau}\in M$ such that:
$$
t_{\tau}\beta b_{\tau} \in \alpha A ~.
$$

Putting:
$$
s:=\prod_{\sigma\in S} s_{\sigma} ~,\qquad t:=\prod_{\tau\in T} t_{\tau} ~,
$$
we conclude that:
$$
s\alpha A \subseteq \beta B ~,\qquad t\beta B \subseteq \alpha A ~,
$$
and so consequently, $xB \subseteq A$ and $yA \subseteq B$ if we put:
$$
x:=t\cdot\frac{\beta}{\alpha} ~,\qquad y:=s\cdot\frac{\alpha}{\beta} ~.
$$

We have then $x\in (A:B)$, $y\in (B:A)$, and $xy = st\in M$ so that $xy$ is a non-negative power of $\lambda$:
$$
xy = \lambda^k
$$
for some $k\in\Z_{\ge 0}$.
\smallskip

If now additionally the hypotheses of {\it (ii)} are fulfilled then the ideals $A$ and $B$ are invertible ideals of $O$. We can then write $(A :B) = A B^{-1}$, and since now
$$
O\cdot x \subseteq A B^{-1}
$$
we have
$$
O\cdot x = A B^{-1}\cdot U
$$
with a certain ideal $U$ of $O$ (namely, $U = A^{-1} B \cdot O x$). Similarly,
$$
O\cdot y = A^{-1} B\cdot V
$$
with a certain ideal $V$ of $O$. Then $U V = O\cdot xy = O\cdot\lambda^k$ which shows first that $U$ and $V$ are both prime to the conductor of $O$ (since $O\cdot\lambda$ is), and then that their prime divisors are all among the prime divisors $\{ P_i\}_{i\in I}$ of $\lambda$ in $O$. Hence,
$$
[U] \in \langle [P_i] \mid ~i\in I\rangle
$$
in $\mathrm{Pic}(O)$; furthermore, as $O\cdot x = A B^{-1}\cdot U$, we have
$$
[A] - [B] + [U] = 0
$$
in $\mathrm{Pic}(O)$. We have shown the `only if' part of {\it (ii)}.
\smallskip

For the proof of {\it (iii)} observe that we clearly have:
$$
xO_K B \subseteq O_K A ~,\qquad yO_K A \subseteq O_K B ~.
$$

Since all (nonzero) ideals of $O_K$ are invertible, we can repeat the above arguments, substituting $A$ and $B$ by $O_K A$ and $O_K B$, respectively, to obtain the conclusion of {\it (iii)}.
\bigskip

Let us then prove the `if' part of {\it (i)}: Suppose that we have elements $x\in (A :B)$ and $y\in (B :A)$ such that $xy=\lambda^k$ for some non-negative integer $k$.

As $x\in (A :B)$ we certainly have $x\B = xO_{(M)}\cdot B \subseteq O_{(M)}\cdot A = \A$. On the other hand, since $y\in (B :A)$ and since $xy = \lambda^k$ is a unit in $O_{(M)}$ we have:
$$
\A = O_{(M)}\cdot A = O_{(M)}\cdot xy\cdot A \subseteq O_{(M)}\cdot x B = x\B ~.
$$

Hence, $\A=x\B$ and $\A\sim \B$.
\bigskip

Finally, we prove the `if' part of {\it (ii)}: By assumption there are then non-negative integers $v_i$ for $i\in I$ such that:
$$
A \sim B \cdot \prod_{i\in I} P_i^{v_i} \leqno{(\ast)}
$$
(we can choose the $v_i$ to be negative since $\mathrm{Pic}(O)$ is finite).

Let
$$
O\cdot\lambda = \prod_{i\in I} P_i^{m_i}
$$
be the prime factorization of $O\cdot\lambda$; by assumption, each $m_i$ is nonzero. Choose then $k\in\N$ such that $k\cdot m_i \ge v_i$ for each $i$, put $u_i := k\cdot m_i - v_i$, and:
$$
U := \prod_{i\in I} P_i^{u_i} ~,\qquad V := \prod_{i\in I} P_i^{v_i} ~;
$$
these are ideals of $O$ prime to the conductor, and we have:
$$
U V = \prod_{i\in I} P_i^{u_i+v_i} = \prod_{i\in I} P_i^{km_i} = O\cdot \lambda^k ~.\leqno{(\ast\ast)}
$$

Now, by $(\ast)$ and the definition of $V$ we have
$$
A^{-1} B V = O\cdot y
$$
for some $y\in K^{\times}$. Since $V$ is an ideal of $O$ we have:
$$
y \in A^{-1} B V \subseteq A^{-1} B = (B:A) ~.
$$

Also, $[A B^{-1}] = [V] = [U^{-1}]$ in $\mathrm{Pic}(O)$ because of $(\ast\ast)$; hence,
$$
A B^{-1} U = O\cdot x
$$
for some $x\in K^{\times}$. Since $U$ is an ideal of $O$ (as all $u_i$ are $\ge 0$), we see that:
$$
x\in A B^{-1} U \subseteq A B^{-1} = (A:B) ~.
$$

Now, $O\cdot xy = U V = O\cdot\lambda^k$ by $(\ast\ast)$; changing $x$ by a unit of $O$ if necessary we then have $xy = \lambda^k$. By the already proved `if' part of {\it (i)} we conclude that $\A\sim \B$.
\end{proof}

\noindent {\bf Remarks:} In the setting of Theorem \ref{t1}, the ideals $A$ and $B$, as well as the fractional ideals $(A:B)$ and $(B:A)$ are all finitely generated abelian groups of rank $[K:\Q]$. If $A$ and $B$ are given explicitly via generators then generators for $(A:B)$ and $(B:A)$ can be computed.
\smallskip

The fractional ideals $(A:B)$ and $(B:A)$ are in particular finitely generated modules over the order $O$, and $O$-module generators can be found if the ideals $A$ and $B$ are given explicitly. Thus, the question of solvability of a single equation $xy=\lambda^k$ with $x\in (A:B)$ and $y\in (B:A)$ reduces to the question of solvability in the order $O$ of a single quadratic equation
$$
f(x_1,\ldots,x_s) = \lambda^k
$$
where $f$ is a quadratic form with coefficients in $K$ that can be determined algorithmically when $A$ and $B$ are explicitly known.
\smallskip

In \cite{gs} it was remarked that the methods of that work show that there is an algorithm for deciding a question like this, i.e., the question of solvability of a quadratic equation in an explicitly given order of an algebraic number field.
\smallskip

Hence, condition $(i)$ of Theorem \ref{t1} would become an algorithmically decidable criterion if one could somehow limit the $k$'s that have to be considered to a finite number. In a sense, such a reduction to consideration of only finitely many $k$'s is what is happening under the favorable conditions of $(ii)$ of the theorem, the main point being the finiteness of $\mathrm{Pic}(O)$.
\smallskip

The question of whether condition $(i)$ is algorithmically decidable in the general case where one or more of the ideals $A$, $B$, and $O\cdot \lambda$ are not prime to the conductor of $O$ is a more complicated question that we will return to elsewhere.

\section{The experiment.}\label{comp}
To give a quantitative description of the explanatory power of our adjusted invariant we investigate it in the context of \cite{lund_jensen} (a Master's Thesis written under the supervision of the first author.) In this work, a large-scale experiment was performed to investigate how close one would get to understanding strong shift equivalence in the set $\uni$ consisting of all irreducible $2\times 2$ matrices $A$ with integer entries and entry sum $\leq 25$, by combining the known invariants with a brute force search for elementary shift equivalences. The invariants available for this project were

\begin{enumerate}
\item The essential Jordan form of $A$, disregarding the null space of $A$, if necessary
\item The Bowen-Franks type groups $\Z^n/p(A)\Z^n$
where $p$ is one of
\begin{gather*}
 x \pm 1,		2x \pm 1,
	x^2 \pm x \pm 1,	
x^2 \pm 2x + 1,	x^2 \pm 1,	
2x^2 \pm x - 1,\\	2x^2 \pm 3x + 1,
	4x^2 \pm 4x + 1,	
4x^2 - 1.
\end{gather*}
\item The BMT invariant under the assumption that $\lambda$ be a unit of the quadratic number field $\Q(\lambda)$.
\end{enumerate}

The condition that $\lambda$ be a unit is of course equivalent to  $\Z[1/\lambda]\subseteq \Z[\lambda]$, and hence in this case the BMT invariant coincides with the ideal class invariant resulting from {\it (iii)} of Theorem \ref{t1}. Obviously, this is a very strong restriction on $\lambda$, and the present work arose initially out of a desire to remove this restriction.
\smallskip

We will say that $A\notknown B$ when all of the above invariants coincide whenever they are defined. As above, when $A$ and $B$ are strong shift equivalent, we write $A\approx B$, and when a concrete strong shift equivalence from $A$ to $B$ is known to us, we write $A\known B$.

We obviously have
\[
A\known B\Longrightarrow
A\approx B\Longrightarrow
A\notknown B,
\]
where the equivalence relation ``$\approx$'' is the one we are interested in, but do not know how to decide. We hence try to approximate the relation by coarser and finer equivalence relations which may be decided, in the case of ``$\known$'' by looking up the pair $(A,B)$ in the database obtained in \cite{lund_jensen} (the database is publicly accessible at http://www.math.ku.dk/symbdyn/), and in the case of $\notknown$ by computing and comparing the invariants. The work in \cite{lund_jensen} resulted in partitions with
$|\uni/\known|=3522$
and
$|\uni/\notknown|=2068$.
It is expected that the rather large gap results from the fact that Lund Jensen did not have the computer resources to perform a complete search for elementary shift equivalences in relevant $3\times 3$-matrices, and indeed a focused search on ``hard cases'' using a variation of Baker's method (\cite{baker2}) might be employed to decrease the upper bound.

The invariant presented in the paper at hand allows us to increase the lower bound by proving --  in the cases where ${\mathbb Z}[1/\lambda]$ differs from ${\mathbb Z}[\lambda]$ to which no invariant was available to Lund Jensen -- that certain pairs of matrices are \textbf{not} strong shift equivalent. A total of $29$ classes of $\uni/\notknown$ could be distinguished this way.

For example, previously we did not know how to tell the following three matrices
\[
\begin{bmatrix}
5&	13\\	6&	1
\end{bmatrix}
\begin{bmatrix}
5&	6\\	13&	1
\end{bmatrix}
\begin{bmatrix}
4&	9\\	9&	2
\end{bmatrix}
\]
apart, but our new invariant proves that none of them are strong shift equivalent to another.
\smallskip

As another concrete example consider the matrices:
$$
A=\begin{bmatrix}14& 2\\ 1 &0\end{bmatrix}\qquad
B=\begin{bmatrix}13&5\\ 3& 1\end{bmatrix}.
$$

In the large scale computation described above, these matrices turned out to have (the same Jordan forms and) the same Bowen-Franks invariant w.r.t.\ the polynomials listed above. Hence, at that point strong equivalence could not be excluded. But the invariant of part {\it (ii)} (and, in this case, equivalently part {\it (iii)}) of Theorem \ref{t1} does show that the matrices are not strong shift equivalent. In the large scale experiment described above we used a Magma (cf.\ \cite{magma}) script to check the condition of part {\it (iii)} of Theorem \ref{t1} in unresolved cases, but for the concrete example at hand, we can give an explicit, manual verification:

The matrices have characteristic polynomial $x^2-14x-2$ with roots $7\pm \sqrt{51}$. Thus, we put $\lambda := 7+\sqrt{51}$ and consider the quadratic field $K=\Q(\sqrt{51})$ of discriminant $4\cdot 51 = 2^2\cdot 3\cdot 17$. The ring of integers of $K$ is $O_K = \Z + \Z \sqrt{51}$. Thus:
$$
O := \Z[\lambda] = O_K
$$
is in fact the maximal order in this case, and {\it (ii)} of Theorem \ref{t1} applies.

Eigenvectors for $A$ and $B$ w.r.t.\ the eigenvalue $\lambda$ are:
$$
\left( \begin{array}{c} \lambda \\ 1 \end{array} \right) \qquad \mbox{and} \qquad \left( \begin{array}{c} -5 \\ 6 - \sqrt{51} \end{array} \right) ~,
$$
respectively. Hence, we need to consider the following ideals of $O_K$:
$$
\A := O_K \cdot \lambda + O_K \cdot 1 = O_K
$$
and
$$
\B := O_K \cdot 5 + O_K \cdot (6 - \sqrt{51}) = O_K \cdot 5 + O_K \cdot (1 - \sqrt{51}) ~,
$$
and the question now simply becomes whether the ideal $\B$ is principal. If not, the matrices $A$ and $B$ are not strong shift equivalent.

Now, the ideal $\B$ is in fact one of the two distinct prime divisors of $5$ in $K$ (cf.\ e.g.\ Thm.\ $25$ of \cite{marcus}.) Hence the norm of $\B$ is $\Nm_{K/\Q}(\B) = 5$ and the question becomes whether there exists a number in $O_K$ with norm $\pm 5$, i.e.\ whether one of the two equations:
$$
x^2 - 51 y^2 = \pm 5 \leqno{(\ast)}
$$
has a solution in integers $x$ and $y$. If not, $\B$ is not principal.

Now, none of the equations $(\ast)$ has in fact a solution in integers. This can be seen explicitly as follows: The fundamental unit $\epsilon$ of $O_K$ is $\epsilon := 50 + 7\sqrt{51}$ as can be ascertained for instance using the continued fractions method (cf.\ e.g.\ \cite{bs}, p.\ 134.) Thus, by \cite{hasse}, pp. 578--579 for instance, the integral non-solvability of $(\ast)$ follows from its non-solvability with $y$ an integer in the range:
$$
- \frac{\epsilon}{\sqrt{51}} < y < \frac{5+\epsilon}{\sqrt{51}} ~.
$$

As $- \frac{\epsilon}{\sqrt{51}} \approx -14.001$ and $\frac{5+\epsilon}{\sqrt{51}} \approx 14.70$ we see that it is enough to verify that none of the numbers $\pm 5 + 51 y^2$ with $y$ running through integers in the interval $1\le y \le 14$ is the square of an integer. And this is easily checked of course.


\end{document}